\documentclass[a4 paper,12pt]{article}
\usepackage[T1]{fontenc}
\usepackage[latin1]{inputenc}
\usepackage{hyperref}
\usepackage[dvips]{graphics}
\usepackage[english]{babel}
\usepackage{amsmath}
\usepackage{amssymb}
\usepackage{amsopn}
\usepackage{amsbsy}
\usepackage{amstext}
\usepackage{latexsym}
\usepackage{amsfonts}
\usepackage{array}
\usepackage{epsfig}
\usepackage{mathrsfs}
\usepackage{color}
\definecolor{vert}{rgb}{0,0.5,0}
\title{An analog of the arithmetic triangle obtained by replacing the products by the least common multiples}
\author{\sc Bakir FARHI}
\date{}

\newtheorem{thm}{Theorem}
\newtheorem{prop}[thm]{Proposition}

\newtheorem{coll}[thm]{Corollary}

\let\epsilon=\varepsilon
\def\lcm{{\rm lcm}}

\def\EMdash{\leavevmode\hbox to 7.5mm{\vrule height .63ex depth -.59ex
    width 5.4mm\hfill}}

\begin{document}
\maketitle ~\vspace{-1.5cm}
\begin{center}
{\tt bakir.farhi@gmail.com}
\end{center}

\section{Introduction}~

The Al-Karaji arithmetic triangle is the triangle consisting of
the binomial coefficients $\binom{n}{k}$ ($n , k \in \mathbb{N} ,
n \geq k$). Precisely, for each $n \in \mathbb{N}$, the
$n$\textsuperscript{th} row of that triangle is:
$$\binom{n}{0} ~~ \binom{n}{1} ~~ \dots ~~ \binom{n}{n} ,$$
where
\begin{equation}\label{eq2}
\binom{n}{k} := \frac{n!}{k! (n - k)!} = \frac{n \times (n - 1)
\times \dots \times (n - k + 1)}{1 \times 2 \times \dots \times k}
\end{equation}
So the beginning of the arithmetic (or binomial) triangle is given
by:
$$
\begin{array}{lllllll}
1 & ~ & ~ & ~ & ~ & ~ & ~ \\
1 & 1 ~ & ~ & ~ & ~ & ~ \\
1 & 2 & 1 & ~ & ~ & ~ & ~ \\
1 & 3 & 3 & 1 & ~ & ~ & ~ \\
1 & 4 & 6 & 4 & 1 & ~ & ~ \\
1 & 5 & 10 & 10 & 5 & 1 & ~ \\
1 & 6 & 15 & 20 & 15 & 6 & 1 \\
\hspace{0.7mm}\vdots & \hspace{0.7mm}\vdots & \hspace{0.7mm}\vdots
& \hspace{0.7mm}\vdots & \hspace{0.7mm}\vdots &
\hspace{0.7mm}\vdots & \hspace{5mm}\ddots
\end{array}
$$
Note that the construction of the triangle rests on the property
that each number of a given row is the sum of the numbers which
are situated just above. Explicitly, we have:
\begin{equation}\label{eq1}
\binom{n}{k} = \binom{n-1}{k-1} + \binom{n-1}{k} ~~~~~~ (\forall k
, n ~\text{such that}~ n \geq k \geq 1)
\end{equation}

Historically, the first mathematician who discovered the binomial
triangle was the pioneer arabic mathematician Al-Karaji (953 -
1029 AD). He drew this triangle until its 12\textsuperscript{th}
row and noted the process of its recursive construction by
pointing out (\ref{eq1}). More interestingly, Al-Karaji discovered
the binomial formula:
\begin{equation}\label{eq3}
(x + y)^n = \sum_{k = 0}^{n} \binom{n}{k} x^k y^{n - k} ~~~~
(\forall n \in \mathbb{N})
\end{equation}
After Al-Karaji, several other mathematicians of the Islamic
civilization reproduced that very important triangle (Al-Khayyam,
Al-Samawal, Al-Tusi, Al-Farisi, Ibn Al-Banna, Ibn Munaim,
Al-Kashi, $\dots$). The same triangle have been discovered again
in China (Yang Hui in the 13\textsuperscript{th} century). In
Europ (16\textsuperscript{th} century), several mathematicians
remarked the importance of Al-Karaji's triangle (Stifel,
Tartaglia, Pascal, $\dots$).

In this paper, we are going to obtain the analog of Al-Karaji's
triangle by substituting in Formula (\ref{eq2}) the products by
the least common multiples. If we use the formula $\binom{n}{k} =
\frac{n!}{k! (n - k)!}$, the lcm-analog of the binomial
coefficient $\binom{n}{k}$ would be:
$$\frac{\lcm(1 , 2 , \dots , n)}{\lcm(1 , 2 , \dots , k) \times \lcm(1 , 2 , \dots , n - k)} .$$
But this analogy is not quite interesting because those last
numbers are not all integers. For example, for $n = 6 , k = 3$, we
have:
$$\frac{\lcm(1 , 2 , \dots , 6)}{\lcm(1 , 2 , 3) \times \lcm(1 , 2 , 3)} = \frac{5}{3} \not\in \mathbb{Z} .$$
In order to obtain an interesting analogy, we will use rather the
formula $\binom{n}{k} = \frac{n \times (n - 1) \times \dots \times
(n - k + 1)}{1 \times 2 \times \dots \times k}$. So, the
lcm-analog of a binomial coefficient $\binom{n}{k}$ which we must
consider is:
\begin{equation}\label{eq10}
\left[\begin{array}{c} n \\ k
\end{array}\right] := \frac{\lcm(n , n - 1 , \dots , n - k + 1)}{\lcm(1 , 2 , \dots , k)}
\end{equation}
(We naturally conventione that $\lcm(\emptyset) =
1$).~\vspace{2mm}

Notice that a table of the numbers $[\begin{subarray}{c} n
\\ k \end{subarray}]$ was already given by A. Murthy (2004) and extended by
E. Deutsch (2006) in the On-Line Encyclopedia of Integer Sequences
(see the sequence \underline{\color{blue}{A093430}} of OEIS).
However, to my knowledge, no property was already proved about
those numbers in comparison with their analog binomial numbers.

\section{Results}~

We begin with the easy result showing that the rational numbers
$[\begin{subarray}{c} n
\\ k \end{subarray}]$, defined by (\ref{eq10}), are all integers.
We have the following:
\begin{prop}\label{prop1}
For all natural numbers $n , k$ such that $n \geq k$, the positive
rational number $[\begin{subarray}{c} n \\ k
\end{subarray}]$ is an integer.
\end{prop}

\noindent{\bf Proof.} Let $n , k$ be natural numbers such that $n
\geq k$. Among the $k$ consecutive integers $n , n - 1 , \dots , n
- k + 1$, one at least is a multiple of $1$, one at least is a
multiple of $2$, $\dots$, and one at least is a multiple of $k$.
This implies that $\lcm(n , n - 1 , \dots , n - k + 1)$ is a
multiple of each of the positive integers $1 , 2 , \dots , k$.
Consequently $\lcm(n , n - 1 , \dots , n - k + 1)$ is a multiple
of $\lcm(1 , 2 , \dots , k)$, which confirms that
$[\begin{subarray}{c} n
\\ k \end{subarray}]$ is an integer. The proposition is
proved.\hfill $\blacksquare$~\vspace{2mm}

\noindent{\bf Definition.} Throughout this paper, we call the
numbers $[\begin{subarray}{c} n \\ k
\end{subarray}]$: ``the $\lcm$-binomial numbers'' and we call
the triangle consisting of them: ``the $\lcm$-binomial
triangle''.~\vspace{2mm}

The beginning of the $\lcm$-binomial triangle is given in the
following:
$$
\begin{array}{lllllll}
1 & ~ & ~ & ~ & ~ & ~ & ~ \\
1 & 1 ~ & ~ & ~ & ~ & ~ \\
1 & 2 & 1 & ~ & ~ & ~ & ~ \\
1 & 3 & 3 & 1 & ~ & ~ & ~ \\
1 & 4 & 6 & \color{vert}{2} & 1 & ~ & ~ \\
1 & 5 & 10 & 10 & 5 & 1 & ~ \\
1 & 6 & 15 & \color{vert}{10} & \color{vert}{5} & \color{vert}{1} & 1 \\
\hspace{0.7mm}\vdots & \hspace{0.7mm}\vdots & \hspace{0.7mm}\vdots
& \hspace{0.7mm}\vdots & \hspace{0.7mm}\vdots &
\hspace{0.7mm}\vdots & \hspace{5mm}\ddots
\end{array}
$$
(Here the colored numbers in green are those that are different
from their analog binomial numbers).~\vspace{1mm}

Now, we are going to establish less obvious results concerning the
$\lcm$-binomial numbers.

\begin{thm}\label{t1}
For all natural numbers $n , k$ such that $n \geq k$, the
$\lcm$-binomial number $[\begin{subarray}{c} n \\ k
\end{subarray}]$ divides the binomial number $\binom{n}{k}$.
\end{thm}

\noindent{\bf Proof.} Actually the theorem can be immediately
showed by using a result of S. Hong and Y. Yang \cite{hy} which
states that for all integers $k , n$ (with $k \geq 0$, $n \geq
1$), the positive integer $g_k(1)$ divides the positive integer
$g_k(n)$, where $g_k$ denotes the Farhi arithmetical
function\footnote{By definition: $g_k(n) := \frac{n (n + 1) \cdots
(n + k)}{\lcm(n , n + 1 , \dots , n + k)} ~~~~ (\forall k , n)$.}
(see Lemma 2.4 of \cite{hy}). But in order to put the reader at
their ease, we give in what follows an independent and complete proof.\\
Let $n , k \in \mathbb{N}$ such that $n \geq 1$ and $n \geq k$.
The statement of the theorem is clearly equivalent to the
following inequalities:
\begin{equation}\label{eq4}
v_p\left(\binom{n}{k}\right) \geq v_p\left(\left[\begin{array}{c}
n \\ k
\end{array}\right]\right) ~~~~ \text{(for all prime number $p$)}
\end{equation}
(where $v_p$ denotes the usual $p$-adic valuation).\\
Let us show (\ref{eq4}) for a given prime number $p$. On the one
hand, we have:
\begin{eqnarray}
v_p\left(\binom{n}{k}\right) & = & v_p\left(\frac{n!}{k! (n -
k)!}\right) \notag \\
& = & v_p(n!) - v_p(k!) - v_p((n - k)!) \notag \\
& = & \sum_{\alpha = 1}^{\infty} \left\lfloor \frac{n}{p^{\alpha}}
\right\rfloor - \sum_{\alpha = 1}^{\infty} \left\lfloor
\frac{k}{p^{\alpha}} \right\rfloor - \sum_{\alpha = 1}^{\infty}
\left\lfloor \frac{n -
k}{p^{\alpha}} \right\rfloor \notag \\
& = & \sum_{\alpha = 1}^{\infty}\left(\left\lfloor
\frac{n}{p^{\alpha}}\right\rfloor -
\left\lfloor\frac{k}{p^{\alpha}}\right\rfloor - \left\lfloor
\frac{n - k}{p^{\alpha}}\right\rfloor\right) \label{eq6}
\end{eqnarray}
(where $\lfloor.\rfloor$ represents the integer part function).\\
It is important to stress that each of the terms
$(\lfloor\frac{n}{p^{\alpha}}\rfloor -
\lfloor\frac{k}{p^{\alpha}}\rfloor - \lfloor\frac{n -
k}{p^{\alpha}}\rfloor)$ ($\alpha \geq 1$), of the last sum, is
nonnegative. indeed, for all positive integer $\alpha$, we have:
$$\left\lfloor\frac{k}{p^{\alpha}}\right\rfloor + \left\lfloor\frac{n - k}{p^{\alpha}}\right\rfloor \leq \frac{k}{p^{\alpha}} +
\frac{n - k}{p^{\alpha}} = \frac{n}{p^{\alpha}} .$$ But since
$\lfloor\frac{k}{p^{\alpha}}\rfloor + \lfloor\frac{n -
k}{p^{\alpha}}\rfloor$ is an integer, then we have even:
$$\left\lfloor\frac{k}{p^{\alpha}}\right\rfloor + \left\lfloor\frac{n - k}{p^{\alpha}}\right\rfloor \leq \left\lfloor\frac{n}{p^{\alpha}}\right\rfloor ,$$
which confirms the stressed fact.\\
Now, on the other hand, we have:
\begin{eqnarray*}
v_p\left(\left[\begin{array}{c} n \\ k \end{array}\right]\right) &
= & v_p\left(\frac{\lcm(n , n - 1 , \dots , n - k + 1)}{\lcm(1 , 2
, \dots , k)}\right) \\
& = & a - b ,
\end{eqnarray*}
where
\begin{eqnarray*}
a & := & v_p(\lcm(n , n - 1 , \dots , n - k + 1)) ~~\text{and} \\
b & := & v_p(\lcm(1 , 2 , \dots , k)) .
\end{eqnarray*}
Note that because $[\begin{subarray}{c} n \\ k\end{subarray}]$ is
an integer (according to Proposition \ref{prop1}), we have $a \geq
b .$ \\
By definition, $a$ is the greatest exponent $\alpha$ of $p$
for which $p^{\alpha}$ divides at least an integer of the range
$(n - k , n]$. Since for all $\alpha \in \mathbb{N}$, the number
of integers belonging to the range $(n - k , n]$, which are
multiples of $p^{\alpha}$, is exactly equal to
$\lfloor\frac{n}{p^{\alpha}}\rfloor - \lfloor\frac{n -
k}{p^{\alpha}}\rfloor$, then we have:
\begin{equation}\label{eq7}
a = \max\left\{\alpha \in \mathbb{N} :
\left\lfloor\frac{n}{p^{\alpha}}\right\rfloor -
\left\lfloor\frac{n - k}{p^{\alpha}}\right\rfloor \geq 1\right\}
\end{equation}
Similarly, $b$ is (by definition) the greatest exponent $\alpha$
of $p$ for which $p^{\alpha}$ divides at least an integer of the
range $[1 , k]$. But since for all $\alpha \in \mathbb{N}$, the
number of integers belonging to the range $[1 , k]$, which are
multiples of $p^{\alpha}$, is exactly equal to
$\lfloor\frac{k}{p^{\alpha}}\rfloor$, then we have:
\begin{equation}\label{eq8}
b = \max\left\{\alpha \in \mathbb{N} :
\left\lfloor\frac{k}{p^{\alpha}}\right\rfloor \geq 1\right\}
\end{equation}
Remarking that the sequence
${\left(\lfloor\frac{n}{p^{\alpha}}\rfloor - \lfloor\frac{n -
k}{p^{\alpha}}\rfloor\right)}_{\alpha \in \mathbb{N}}$ is
non-increasing (since each of the terms
$\lfloor\frac{n}{p^{\alpha}}\rfloor - \lfloor\frac{n -
k}{p^{\alpha}}\rfloor$ represents the number of integers lying in
the range $(n - k , n]$, which are multiples of $p^{\alpha}$), we
have:
$$\forall \alpha \in \mathbb{N} , \alpha \leq a :~~ \left\lfloor\frac{n}{p^{\alpha}}\right\rfloor - \left\lfloor\frac{n -
k}{p^{\alpha}}\right\rfloor \geq 1 .$$ Further, from the
definition of $b$, we have:
$$\forall \alpha \in \mathbb{N} , \alpha > b :~~ \left\lfloor\frac{k}{p^{\alpha}}\right\rfloor = 0 .$$
Consequently, we have:
$$\forall \alpha \in \mathbb{N} \cap (b , a] :~~ \left\lfloor\frac{n}{p^{\alpha}}\right\rfloor - \left\lfloor\frac{n -
k}{p^{\alpha}}\right\rfloor -
\left\lfloor\frac{k}{p^{\alpha}}\right\rfloor \geq 1 .$$ According
to (\ref{eq6}), it follows that:
\begin{eqnarray*}
v_p\left(\binom{n}{k}\right) & = & \sum_{\alpha =
1}^{\infty}\left(\left\lfloor \frac{n}{p^{\alpha}}\right\rfloor -
\left\lfloor\frac{n - k}{p^{\alpha}}\right\rfloor - \left\lfloor
\frac{k}{p^{\alpha}}\right\rfloor\right) \\
& \geq & \sum_{b < \alpha \leq a}\left(\left\lfloor
\frac{n}{p^{\alpha}}\right\rfloor - \left\lfloor\frac{n -
k}{p^{\alpha}}\right\rfloor - \left\lfloor
\frac{k}{p^{\alpha}}\right\rfloor\right) \\
& \geq & \sum_{b < \alpha \leq a} 1 \\
& = & a - b \\
& = & v_p\left(\left[\begin{array}{c} n
\\ k\end{array}\right]\right) ,
\end{eqnarray*}
which confirms (\ref{eq4}) and completes this
proof.\hfill$\blacksquare$~\vspace{2mm}

Now, by Theorem \ref{t1}, we see that the ratios
$\binom{n}{k}/[\begin{subarray}{c} n \\k \end{subarray}]$ are
actually positive integers. But it certainly remains several other
profound properties to discover about those numbers. We can ask
for example about the couples $(n , k)$ satisfying the equality
$\binom{n}{k} =
[\begin{subarray}{c} n \\ k \end{subarray}]$.\\
The following theorem shows a very important property for the
ratios $\binom{n}{k}/[\begin{subarray}{c} n \\k \end{subarray}]$.
We derive from it for example that for a fixed column $k$, the
numbers ${(\binom{n}{k}/[\begin{subarray}{c} n \\k
\end{subarray}])}_{n \geq k}$ lie in a finite set of positive
integers.
\begin{thm}\label{t2}
For all $k \in \mathbb{N}$, the sequence of positive integers
${\left(\frac{\binom{n}{k}}{\left[\begin{subarray}{c} n \\ k
\end{subarray}\right]}\right)}_{\!\!n \geq k}\!\!$ is periodic and its
smallest period $T_k$ is given by:
$$T_k = \prod_{p ~\text{prime, $p < k$}} p^{\alpha_p} ,$$
where
$$\alpha_p = \begin{cases} 0 & \text{if}~ v_p(k) \geq \displaystyle\max_{1 \leq i < k}
v_p(i) \\
\displaystyle\max_{1 \leq i < k} v_p(i) & \text{otherwise}
\end{cases} ~~~~~~~~~~~~~~ (\forall p ~\text{prime}, p < k) .$$
\end{thm}
As an important consequence, we derive the following:
\begin{coll}\label{c1}
For all $k \in \mathbb{N}$, the positive integer $\lcm(1 , 2 ,
\dots , k - 1)$ is a period of the sequence
${\left(\frac{\binom{n}{k}}{\left[\begin{subarray}{c} n \\ k
\end{subarray}\right]}\right)}_{\!\!n \geq k}$.
\end{coll}

Admitting Theorem \ref{t2}, the proof of Corollary \ref{c1}
becomes obvious: it suffices to remark that the exact period
$T_k$, given by Theorem \ref{t2}, of the sequence
${\left(\binom{n}{k}/\left[\begin{subarray}{c} n
\\ k \end{subarray}\right]\right)}_{n \geq k}$ clearly divides $\lcm(1 , 2 , \dots , k -
1)$.

To prove Theorem \ref{t2}, we use the arithmetical functions $g_k$
($k \in \mathbb{N}$) introduced by the author in \cite{f} and
studied later by Hong and Yang \cite{hy} and by Farhi and Kane
\cite{fk}. For a given $k \in \mathbb{N}$, the function $g_k$ is
defined by:
$$
\begin{array}{rll}
g_k : \mathbb{N} \setminus \{0\} & \longrightarrow & \mathbb{N}
\setminus \{0\} \\
n & \longmapsto & g_k(n) := \frac{n (n + 1) \cdots (n + k)}{\lcm(n
, n + 1 , \dots , n + k)} ~\cdot
\end{array}
$$
In \cite{f}, it is just remarked that $g_k$ is periodic and that
$k!$ is a period of $g_k$. Then Hong and Yang \cite{hy} improved
that period to $\lcm(1 , 2 , \dots , k)$ and recently, Farhi and
Kane \cite{fk} have obtained the exact period of $g_k$ which is
given by:
$$
P_k = \prod_{p \ \textrm{prime}, \ p\leq k} p^{\begin{cases}0 &
\textrm{if} \ v_p(k+1)\geq \max_{1\leq i\leq k} v_p(i)\\
\max_{1\leq i\leq k} v_p(i) & \textrm{otherwise} \end{cases}}.
$$
Knowing this result, the proof of Theorem \ref{t2} becomes
easy:~\vspace{2mm}

\noindent{\bf Proof of Theorem \ref{t2}.} For a fixed $k \in
\mathbb{N}$, a simple calculus shows that for any $n \in
\mathbb{N}$, we have:
$$
\frac{\binom{n}{k}}{\left[\begin{subarray}{c} n \\
k\end{subarray}\right]} = \frac{g_{k - 1}(n - k + 1)}{g_{k -
1}(1)} .
$$
This last identity clearly shows that for any given $k \in
\mathbb{N}$, the sequence
${\left(\binom{n}{k}/[\begin{subarray}{c} n
\\ k
\end{subarray}]\right)}_{n \geq k}$ is periodic and that its exact
period is equal to the exact period of $g_{k - 1}$. So by the
Farhi-Kane theorem, the exact period of
${\left(\binom{n}{k}/[\begin{subarray}{c} n
\\ k
\end{subarray}]\right)}_{n \geq k}$ is $P_{k - 1}$, as claimed in
Theorem \ref{t2}.\hfill$\blacksquare$~\vspace{2mm}

We end this section by giving the $\lcm$-binomial triangle until
its 12\textsuperscript{th} row.
$$
\begin{array}{lllllllllllll}
1 & ~ & ~ & ~ & ~ & ~ & ~ & ~ & ~ & ~ & ~ & ~ & ~ \\
1 & 1 ~ & ~ & ~ & ~ & ~ & ~ & ~ & ~ & ~ & ~ & ~ \\
1 & 2 & 1 & ~ & ~ & ~ & ~ & ~ & ~ & ~ & ~ & ~ & ~ \\
1 & 3 & 3 & 1 & ~ & ~ & ~ & ~ & ~ & ~ & ~ & ~ & ~ \\
1 & 4 & 6 & \color{vert}{2} & 1 & ~ & ~ & ~ & ~ & ~ & ~ & ~ & ~ \\
1 & 5 & 10 & 10 & 5 & 1 & ~ & ~ & ~ & ~ & ~ & ~ & ~\\
1 & 6 & 15 & \color{vert}{10} & \color{vert}{5} & \color{vert}{1} & 1 & ~ & ~ & ~ & ~ & ~ & ~ \\
1 & 7 & 21 & 35 & 35 & \color{vert}{7} & 7 & 1 & ~ & ~ & ~ & ~ &
~ \\
1 & 8 & 28 & \color{vert}{28} & 70 & \color{vert}{14} &
\color{vert}{14} & \color{vert}{2} & 1 & ~ & ~ & ~ & ~ \\
1 & 9 & 36 & 84 & \color{vert}{42} & \color{vert}{42} &
\color{vert}{42} & \color{vert}{6} & \color{vert}{3} & 1 & ~ &
~ & ~ \\
1 & 10 & 45 & \color{vert}{60} & 210 & \color{vert}{42} &
\color{vert}{42} & \color{vert}{6} & \color{vert}{3} &
\color{vert}{1} & 1 & ~ & ~ \\
1 & 11 & 55 & 165 & 330 & 462 & 462 & \color{vert}{66} &
\color{vert}{33} & \color{vert}{11} & 11 & 1 & ~ \\
1 & 12 & 66 & \color{vert}{110} & \color{vert}{165} &
\color{vert}{66} & \color{vert}{462} & \color{vert}{66} &
\color{vert}{33} & \color{vert}{11} & \color{vert}{11} &
\color{vert}{1} & 1 \\
\hspace{0.7mm}\vdots & \hspace{0.7mm}\vdots & \hspace{0.7mm}\vdots
& \hspace{0.7mm}\vdots & \hspace{0.7mm}\vdots &
\hspace{0.7mm}\vdots & \hspace{0.7mm}\vdots & \hspace{0.7mm}\vdots
& \hspace{0.7mm}\vdots & \hspace{0.7mm}\vdots &
\hspace{0.7mm}\vdots & \hspace{0.7mm}\vdots & \hspace{5mm}\ddots
\end{array}
$$
$$
\text{\color{red}{The $\lcm$-analog of Al-Karaji's triangle}}
$$
Note that The $\lcm$-binomial numbers colored in green are those
that are different from their analog binomial numbers.

\section{Some remarks and open problems about the $\lcm$-binomial
numbers}~
\begin{description}
\item[1)] Can we prove Theorem \ref{t1} without use prime number
arguments?
 \item[2)] Describe the set of all the couples $(n , k)$ ($n \geq k \geq
 0$) satisfying $[\begin{subarray}{c} n \\ k \end{subarray}] =
 \binom{n}{k}$.
\item[3)] Let $n \in \mathbb{N}$. Since for any $k \in \{0 , 1 ,
 \dots , n\}$, we have $[\begin{subarray}{c} n \\ k \end{subarray}] \leq
 \binom{n}{k}$ (because $[\begin{subarray}{c} n \\ k \end{subarray}]$ divides $\binom{n}{k}$, according to Theorem
 \ref{t1}) then for all nonnegative real number $x$, we have:
$$
\sum_{k = 0}^{n} \left[\begin{array}{c} n \\ k \end{array}\right]
x^k ~\leq~ \sum_{k = 0}^{n} \binom{n}{k} x^k ~=~ (1 + x)^n ,
$$
that is:
\begin{equation}\label{eq9}
\sum_{k = 0}^{n} \left[\begin{array}{c} n \\ k
\end{array}\right] x^k ~\leq~ (1 + x)^n ~~~~~~ (\forall x \geq 0)~.
\end{equation}
Taking $x = 1$ in (\ref{eq9}), we deduce in particular that for
all $n \in \mathbb{N}$, we have $[\begin{subarray}{c} n
\\ \lceil n/2 \rceil \end{subarray}] \leq 2^n$ (where $\lceil . \rceil$ denotes the ceiling function). But since
$[\begin{subarray}{c} n \\ \lceil n/2 \rceil
\end{subarray}] = \frac{\lcm(n , n - 1 , \dots , n - \lceil n/2 \rceil + 1)}{\lcm(1 , 2 , \dots , \lceil n/2
\rceil)}$ is an integer (according to Proposition \ref{prop1}),
then $\lcm(n , n - 1 , \dots , n - \lceil n/2 \rceil + 1)$ is a
multiple of $\lcm(1 , 2 , \dots , \lceil n/2 \rceil)$.
Consequently we have $\lcm(n , n - 1 , \dots , n - \lceil n/2
\rceil + 1) = \lcm(n , n - 1 ,$ $ \dots , n - \lceil n/2 \rceil +
1; 1 , 2 , \dots , \lceil n/2 \rceil) = \lcm(1 , 2 , \dots , n)$.
So $[\begin{subarray}{c} n \\ \lceil n/2 \rceil \end{subarray}]
\leq 2^n$ gives:
$$\lcm(1 , 2 , \dots , n) \leq 2^n \lcm(1 , 2 , \dots \lceil n/2 \rceil) ~~~~ (\forall n \in \mathbb{N}) .$$
The iteration of the last inequality gives:
$$\lcm(1 , 2 , \dots , n) \leq 2^{n + \lceil n/2 \rceil + \lceil n/4 \rceil + \dots} \leq 2^{2 n + \log_2(n)} = n 4^n ~~~~(\forall n \geq 1) .$$
Hence:
$$\lcm(1 , 2 , \dots , n) \leq n 4^n ~~~~ (\forall n \geq 1) ,$$
which is a nontrivial upper bound of $\lcm(1 , 2 , \dots , n)$.\\
The question which we pose is the following:
\begin{quote}
Can we more judiciously use Relation (\ref{eq9}) to prove a
nontrivial upper bound for the least common multiple of
consecutive integers that is significatively better than the
previous one?
\end{quote}
\item[4)] It is easy to see that unfortunately there is no an
internal composition law $\star$ of $\mathbb{N}$ which satisfies
for any positive integers $n , k$ ($n \geq k$):
$$
\left[\begin{array}{c} n \\ k \end{array}\right] ~=~
\left[\begin{array}{c} n - 1 \\ k - 1 \end{array}\right] \star
\left[\begin{array}{c} n - 1 \\ k \end{array}\right]
$$
(the analog of (\ref{eq1})).\\
Indeed, if we suppose that such a law $\star$ exists then we would
have on the one hand $[\begin{subarray}{c} 2
\\ 1
\end{subarray}] \star [\begin{subarray}{c} 2 \\ 2 \end{subarray}]
= [\begin{subarray}{c} 3 \\ 2 \end{subarray}]$, that is $2 \star 1
= 3$ and on the other hand $[\begin{subarray}{c} 4 \\ 3
\end{subarray}] \star [\begin{subarray}{c} 4 \\ 4 \end{subarray}] = [\begin{subarray}{c} 5 \\ 4
\end{subarray}]$, that is $2 \star 1 = 5$; which gives a
contradiction.\\
The problem which we pose is the following:
\begin{quote}
Find an iterative construction (i.e., a construction row by row)
for the $\lcm$-binomial triangle.
\end{quote}
 \item[5)] For a given positive integer $d$, let
 $\Omega(d)$ denote the number of prime factors of $d$, counting with
 their multiplicities.\\
In this item, we look at the diagonals of the
 $\lcm$-binomial triangle. We constat that the first diagonal (which we note by
 $D_0$) contains only the $1$'s; in other words, we have:
 $$\forall d \in D_0 : \Omega(d) = 0 \leq 0 .$$
 The second diagonal (noted $D_1$) is consisted only on the $1$'s
 and the prime numbers; in other words, we have:
 $$\forall d \in D_1 : \Omega(d) \leq 1 .$$
 Also, the third diagonal of the $\lcm$-binomial triangle (noted $D_2$) is consisted of positive integers having at most two prime factors (counting
 with their multiplicities); in other words, we have:
 $$\forall d \in D_2 : \Omega(d) \leq 2 .$$
 More generally, we have the following:
 \begin{prop}\label{prop2}
For $k \in \mathbb{N}$, let $D_k$ denote the $(k +
1)$\textsuperscript{th} diagonal of the $\lcm$-binomial triangle.
Then, we have:
$$\forall d \in D_k :~~~~ \Omega(d) \leq k .$$
 \end{prop}
 The proof of this proposition is actually very easy and leans
 only on the following simple fact:
 $$\forall n \in \mathbb{N} : \frac{\lcm(1 , 2 , \dots , n , n + 1)}{\lcm(1 , 2 , \dots , n)} = \begin{cases}
 p & \text{if $n + 1$ is a power of a prime $p$} \\
1 & \text{otherwise}
 \end{cases} .$$

 \noindent{\bf Proof of Proposition \ref{prop2}.} Let $k \in
 \mathbb{N}$ fixed and let $d \in D_k$. So, we can write $d$ as:
 $d = [\begin{subarray}{c} n + k \\ n \end{subarray}] = \frac{\lcm(k + 1 , k + 2 , \dots , k + n)}{\lcm(1 , 2 , \dots , n)}$ (for some $n \in
 \mathbb{N}$). It follows that $d$ divides the positive integer $\frac{\lcm(1 , 2 , \dots , n + k)}{\lcm(1 , 2 , \dots ,
 n)}$. But we constat that the last number is the product of the
 $k$ positive integers $\frac{\lcm(1 , 2 , \dots , n + i)}{\lcm(1 , 2 , \dots , n +
 i - 1)}$ $(1 \leq i \leq k)$ each of which is either a
 prime number or equal to $1$ (according to the fact mentioned just before this proof). So,
 it follows that:
 $$\Omega(d) \leq \Omega\left(\frac{\lcm(1 , 2 , \dots , n + k)}{\lcm(1 , 2 , \dots , n)}\right) \leq k .$$
 The proposition is proved.\hfill$\blacksquare$~\vspace{2mm}
\end{description}

Note that by using prime number theory, we can improve the obvious
upper bound of Proposition \ref{prop2} to:
$$\forall d \in D_k : ~~~~ \Omega(d) \leq \mathrm{c}\frac{k}{\log{k}} ,$$
where $\mathrm{c}$ is an absolute positive constant (effectively
calculable).


\begin{thebibliography}{00}

\bibitem{f}
{\sc B. Farhi.} Nontrivial lower bounds for the least common
multiple of some finite sequences of integers, {\it J. Number
Theory}, {\bf 125} (2007), p. 393-411.
\bibitem{fk}
{\sc B. Farhi \& D. Kane.} New results on the least common
multiple of consecutive integers, {\it Proc. Am. Math. Soc}, {\bf
137} (2009), p. 1933-1939.
\bibitem{hy}
{\sc S. Hong \& Y. Yang.} On the periodicity of an arithmetical
function, {\it C. R. Acad. Sci. Paris}, S\'er. {\bf I 346} (2008),
p. 717-721.
\end{thebibliography}
\end{document}